\documentclass{l4dc2027}
\usepackage{times}
\usepackage{mathtools}
\usepackage{graphicx}
\usepackage{placeins}
\usepackage{comment}
\title[Set-based State Estimation for Unknown Dynamical Systems]{
Set-based State Estimation for Unknown Dynamical Systems from Noisy Output Data with Probabilistic Guarantees}
\coltauthor{\Name{Georgios Pantazis}\and\Name{Michelle S. Chong}\\
\addr Department of Mechanical Engineering, Eindhoven University of Technology}
\newtheorem{assumption}{Assumption}
\newcommand{\R}{\mathbb R}
\newcommand{\Pp}{\mathbb P}
\newcommand{\XX}{\mathcal X}
\newcommand{\UU}{\mathcal U}

\newcommand{\T}{\mathcal T}
\newcommand{\Yh}{\widehat{\mathcal Y}}
\newcommand{\Xh}{\widehat{\mathcal X}}
\newcommand{\Vk}{\mathcal V_k}

\usepackage{xcolor}
\usepackage{empheq}

\DeclareMathOperator{\diam}{diam}

\DeclareMathOperator{\rank}{rank}

\begin{document}
\maketitle
\begin{abstract}
We develop a methodology for set-based state
estimation from noisy input-output trajectories,
without a predefined model structure for the
unknown state dynamics or access to training state measurements. We derive a  bound on the state set
diameter that depends on the sensor map and leverage it  as a surrogate to formulate a tractable linear program balancing state set tightness against output trajectory coverage.
 Combining this certificate with deterministic or
probabilistically certified measurement noise sets and leveraging the results from scenario optimization with relaxation, we provide joint a posteriori probabilistic guarantees for state inclusion and state set diameter. Numerical simulations corroborate our developed theoretical guarantees.
\end{abstract}

\section{Introduction}
\label{sec:introduction}
State estimation aims to reconstruct the unmeasured 
states of dynamical systems from available input
and output measurements. Conventional approaches rely on model knowledge and include observer design \cite{Luenberger1971,BernardAndrieuAstolfi2022},
probabilistic filtering \cite{Arulampalam2002,SarkkaSvensson2023},
and set-based estimation, which characterizes
uncertainty through admissible state sets \cite{BertsekasRhodes1971,MilaneseVicino1991}. 
Despite advances in zonotopic set representations
and stochastic/set-membership filtering
\cite{ScottEtAl2016,RegoEtAl2021,
dePaulaEtAl2022,Combastel2016},
a fundamental challenge remains: ensuring
inclusion of the true state, either deterministically
under given uncertainty bounds or with high
probability if data-driven methods are employed,
while maintaining sufficiently tight estimated sets.
Achieving both objectives is particularly
challenging when the dynamics are unknown
and only noisy input-output data are available.

A promising method for set-based estimation with
unknown dynamics is data-driven reachability analysis,
which constructs reachable set over-approximations
in the state space from input-state data
\cite{AlanwarKoch2021,zhang2025data}
or performs set-based state estimation from noisy
input-output data \cite{AlanwarEtAl2022}. However, these guarantees require structural
assumptions on the dynamics, such as linearity
or predefined nonlinear terms, together
with known noise bounds. Such reachability analysis
methods additionally require training using full state measurements. Probabilistic methods quantify the reliability of
data-driven reachable sets through scenario
optimization \cite{DevonportArcak2020,DietrichEtAl2024,rickard2026data,pantazis2026},
Christoffel functions \cite{DevonportEtAl2023},
and ellipsoidal approximations \cite{ParkEtAl2024}.
Other approaches provide probabilistic guarantees
based on the holdout method
\cite{DietrichEtAl2025} or conformal inference
\cite{TebjouEtAl2023,HashemiEtAl2023,NathEtAl2026}.
 Without a predefined model structure for the
unknown state dynamics, or training state measurements, constructing tight state sets only from noisy
input-output trajectories remains challenging. Specifically, the states remain unmeasured during training, and output coverage alone does
not certify both state inclusion and a tight state set size.
In this paper, we address this challenge assuming a known linear sensing map and either
deterministic noise bounds or a
specified noise distribution with possibly unbounded support.  Our contributions are as follows:
\begin{enumerate}
 \item We develop a methodology that learns rectangular
output sets as functions of the input history
from sampled input-output trajectories.
We derive a bound on the diameter of the
corresponding state sets and use its
parameter-dependent term as a surrogate in a
linear program designed to learn the output sets, which directly balances state set tightness against output
trajectory coverage.
We then construct the state sets by taking the
inverse image of the learned output sets, simultaneously accounting for deterministic or probabilistic measurement noise bounds.
Unlike \citet{pantazis2026}, the formulation
requires no state trajectories for training.
\item We establish a posteriori
probabilistic guarantees that the estimated state
sets simultaneously contain the true states and
satisfy the learned diameter bounds over the
entire horizon.
Specifically, since output inclusion alone does not ensure
satisfaction of the learned diameter thresholds,
we include both conditions in the appropriateness
criterion and establish the consistency required
to apply the results of \citet{GarattiCampi2025},
thus avoiding the need to impose a nondegeneracy
assumption.
The resulting distribution-free certificate
depends on the complexity of the learned solution.
Combining this certificate with deterministic
or probabilistically certified measurement noise
sets, we provide explicit bounds on the probability of violating either state inclusion or the state set diameter bounds.
Learning and certification use the same input-output
trajectories, without the need of a separate calibration set.
\end{enumerate}

\emph{Notation}:  $A\oplus B:=\{a+b:a\in A,b\in B\}$, $-B:=\{-b:b\in B\}$, and $h^{-1}(B):=\{x\in\XX:h(x)\in B\}$, $\diam(B):=\sup_{x,z\in B}\|x-z\|_2$ for nonempty $B$. The nullspace of $C$ is $\ker C:=\{q\in\mathbb R^{n_x}:Cq=0\}$. For $i=1,\ldots,n_y$, let $e_i\in\mathbb R^{n_y}$
denote the $i$-th standard basis vector, whose $i$-th
entry is one and all other entries are zero.  Let $\mathcal{T}:=\{0,\ldots,T\}$ and
$\mathcal{T}_{-}:=\{0,\ldots,T-1\}$ denote the time horizon
index sets and consider the sensor index set $\mathcal S:=\{1,\ldots,S\}$. The operator $\operatorname{col}$ denotes the vertical concatenation of elements in a vector.
For a sequence $z(k)\in\mathbb R^d$ and integers
$a\leq b$, define
$z_{a:b}:=\operatorname{col}(z(a),\ldots,z(b))
\in\mathbb R^{(b-a+1)d}$,
where the indices run in increasing order from $a$ to $b$. We finally define the nonnegative integers as $\mathbb Z_{\geq0}$. The indicator $\mathbf 1_A$ equals one if condition
$A$ holds and zero otherwise.
\section{Problem Formulation}
 Consider the system with the following state and output dynamics:

\begin{subequations}\label{eq:system}
\begin{empheq}[left=\empheqlbrace]{align}
&x(k+1) = f(x(k),u(k),w(k)),
&& k \in \mathcal{T}_{-} ,
\label{eq:state_dynamics}\\
&y^{(s)}(k) = h_s(x(k))+v^{(s)}(k),
&& s \in \mathcal{S}, \ k \in \mathcal{T} .
\label{eq:output_dynamics}
\end{empheq}
\end{subequations}

The state $x(k)\in\mathcal X\subseteq\mathbb R^{n_x}$
evolves according to the unknown dynamics
$f:\mathcal X\times\mathcal U\times\mathcal W
\to\mathcal X$, where
$u(k)\in\mathcal U\subseteq\mathbb R^{n_u}$
is the known input and
$w(k)\in\mathcal W\subseteq\mathbb R^{n_w}$
is the unknown process disturbance.
In our setting, neither the state trajectory $x_{0:T}$ nor
the disturbance $w_{0:T-1}$
 are accessible. The state domain $\mathcal{X}$ is known, and can possibly be the entire
$\mathbb{R}^{n_x}$, whereas $\mathcal{W}$ need not be known. Note that
neither set is assumed to be bounded.
For each sensor $s\in\mathcal S$,
the sensor map is
$h_s:\mathcal X\to\mathbb R^{n_s}$,
and the measured output
$y^{(s)}(k)\in\mathbb R^{n_s}$ satisfies \eqref{eq:output_dynamics}
where $v^{(s)}(k)\in\mathbb R^{n_s}$
denotes the measurement noise for sensor $s$.  For each $k\in\mathcal T$, we thus define the concatenated
output, measurement noise, and sensing map by $y(k):=\operatorname{col}
(y^{(1)}(k),\ldots,y^{(S)}(k))$,
$v(k):=\operatorname{col}
(v^{(1)}(k),\ldots,v^{(S)}(k))$,
$h(x):=\operatorname{col}
(h_1(x),\ldots,h_S(x))$, respectively, 
where $n_y:=\sum_{s=1}^{S}n_s$.
Thus, we have that $y(k),v(k)\in\mathbb R^{n_y}$ and
$h:\mathcal X\to\mathbb R^{n_y}$, and the sensor dynamics can be written in the compact form $y(k)=h(x(k))+v(k)$.   
We now impose the following assumption: 
\begin{assumption} \label{assum:sensor_maps}
For every sensor $s\in \mathcal{S}$, the sensor map is
given by the linear transformation
$h_s(x)=C_sx$, where the matrix $C_s\in\mathbb R^{n_s\times n_x}$ is known. \hfill $\square$
\end{assumption}
Let $C:=\operatorname{col}(C_1,\ldots,C_S)
\in\mathbb R^{n_y\times n_x}$ denote the concatenated sensor matrix, so that $h(x)=Cx$. Writing $u:=u_{0:T-1}$ and $y:=y_{0:T}$, we collect
the input and corresponding output trajectory into a single
vector $\zeta=(u,y)$, which we refer to as a trajectory sample
or simply sample. A collection of such samples is called a multi-sample. We now impose the following assumption:
\begin{assumption} \label{ass:law}
The trajectories in the multi-sample
$D_N=(\zeta^{(1)},\ldots,\zeta^{(N)})$
are independent and identically distributed according
to a possibly unknown probability distribution
$\mathbb P_0$. A new trajectory pair $\zeta=(u, y)$ is also independent from the multi-sample $D_N$ and follows the same probability distribution $\mathbb{P}_0$, i.e., it holds that
$\bigl(\zeta^{(1)},\ldots,\zeta^{(N)},\zeta \bigr)
\sim \prod_{i=1}^{N+1} \mathbb{P}_0.$ \hfill $\square$
\end{assumption}

Assumption~2 permits dependence across time, sensors,
inputs, and outputs within each trajectory.
We further impose the following measurability assumption.

\begin{assumption}\label{assum:maps}
The sets $\mathcal X$, $\mathcal U$, and
$\mathcal W$ are Borel subsets of their
respective Euclidean spaces.
The dynamics
$f:\mathcal X\times\mathcal U\times\mathcal W
\rightarrow\mathcal X$
and sensor maps
$h_s:\mathcal X\rightarrow\mathbb R^{n_s}$,
$s\in\mathcal S$, are Borel measurable.
Furthermore, we consider Borel measurable feedback control
policies of the form $u(k)=\pi_k(y_{0:k},u_{0:k-1})$, with 
$k\in\mathcal{T}_{-}$ and $u_{0:-1}=\emptyset$. \hfill $\square$
\end{assumption}

Assumption~\ref{ass:law} holds for both open-loop
and feedback experiments under identical training
and deployment conditions. Let
$E=(x_0,w_{0:T-1},v_{0:T})\sim\mathbb P_E$ and 
suppose the dynamics, sensor maps, horizon,
uncertainty distribution, and input mechanism
are fixed across experiments. 
Under Assumption~\ref{assum:maps}, either mechanism induces
a measurable trajectory map $\Phi$ such that $\zeta=\Phi(E)$.
Hence, i.i.d.\ uncertainty realizations induce i.i.d.\
input-output trajectories with distribution
$\mathbb P_0(A)=\mathbb P_E(\Phi^{-1}(A))$ for every Borel
set $A$ in the input-output trajectory space.

For the measurement noise, we assume either known deterministic
bounds or a specified distribution $\Pp_v$ that the
noise trajectory $v_{0:T}$ follows. Note that $\Pp_v$ is the marginal
distribution of $v_{0:T}$ under the probability distribution $\Pp_E$.
In the latter case, we generate $N_v$ independent \emph{synthetic}
trajectories from $\Pp_v$, collected in the multi-sample
$D_v=(v_{0:T}^{(1)},\ldots,v_{0:T}^{(N_v)})$,
to construct the noise sets $\Vk(D_v)$.
The distribution of $D_v$ is given by the cartesian product
$\mathbb P_{D_v}:=\prod_{i=1}^{N_v}\mathbb P_v$.
Thus, separate noise measurements during system
operation are not required. Let us now define the probability of exclusion for a new yet unseen noise measurement $\mathbb V_v(D_v)
:=
\Pp_v\left\{
v_{0:T}:\exists k\in\T,\,
v(k)\notin\Vk(D_v)
\right\}$.
We then impose the following assumption:
\begin{assumption}\label{assum:noise_set}
For any $k\in\T$, the mapping
$D_v\mapsto\Vk(D_v)$ is measurable
and the returned set $\Vk(D_v)$ is nonempty and compact.
For a confidence parameter $\beta_v\in[0,1)$,
there exists a violation level
$\epsilon_v(D_v,\beta_v)\in[0,1]$ such that $\Pp_{D_v}\left\{
D_v:\mathbb V_v(D_v)
\leq\epsilon_v(D_v,\beta_v)
\right\}\geq1-\beta_v$.
\hfill $\square$
\end{assumption}
 For known deterministic noise bounds, it holds by construction that
$D_v=\emptyset$ and $\epsilon_v=\beta_v=0$.
For example, using the synthetic trajectories $v^{(j)}_{0:T}$, let
$M:=\max_{j}\max_{k}\|v^{(j)}(k)\|_\infty$ and
$\Vk(D_v):=\{v\in\R^{n_y}:\|v\|_\infty\leq M\}$.
If $Z:=\max_{k\in\T}\|v(k)\|_\infty$ has a continuous
distribution $F_Z$, then it holds that
$\mathbb V_v(D_v)=1-F_Z(M)$.
We then define
$Z_j:=\max_{k\in\mathcal T}\|v^{(j)}(k)\|_\infty$
such that $M=\max_{1\leq j\leq N_v}Z_j$.
For $\epsilon\in(0,1)$, define the quantile
$q_\epsilon:=\inf\{z\in\mathbb R:F_Z(z)\geq1-\epsilon\}$.
Continuity of $F_Z$ then gives
$\mathbb P\{Z_j<q_\epsilon\}=1-\epsilon$.
Moreover, $\mathbb V_v(D_v)>\epsilon$ is equivalent
to $M<q_\epsilon$, meaning that every sampled maximum
$Z_j$ lies below $q_\epsilon$.
Assuming independence of the events,
\begingroup
\setlength{\abovedisplayskip}{1pt}
\setlength{\belowdisplayskip}{1pt}
\setlength{\abovedisplayshortskip}{1pt}
\setlength{\belowdisplayshortskip}{1pt}
\begin{equation*}
\mathbb P_{D_v}\{\mathbb V_v(D_v)>\epsilon\}
=\mathbb P_{D_v}\{Z_j<q_\epsilon
  \text{ for all }j=1,\ldots,N_v\}
=\prod_{j=1}^{N_v}\mathbb P\{Z_j<q_\epsilon\}
=(1-\epsilon)^{N_v}.
\end{equation*}
\endgroup
Choosing $\epsilon_v(D_v,\beta_v)=1-\beta_v^{1/N_v}$
gives $\mathbb P_{D_v}\{\mathbb V_v(D_v)\leq
\epsilon_v(D_v,\beta_v)\}
=1-(1-\epsilon_v(D_v,\beta_v))^{N_v}
=1-\beta_v$,
which establishes Assumption~4.
\section{Characterizing the estimated set of states under unknown dynamics}\label{sec:construction}
For each time $k\in\mathcal T$, let
$\widehat{\mathcal Y}_k(u;\theta)
\subseteq\mathbb R^{n_y}$
denote an estimated output set depending
on the input history $u_{0:k-1}$ and parameterized by the parameters $\theta \in \Theta$. Using the measurement noise sets $\mathcal V_k(D_v)$
introduced in Section~2, we construct the estimated state set:
\begin{equation}\label{eq:preimage}
\begin{aligned}
\Xh_k(u;\theta,D_v)
&:=h^{-1}\bigl(
\Yh_k(u;\theta)\oplus(-\Vk(D_v))
\bigr)=\{x\in\XX:\exists v \in\Vk(D_v),\
Cx+v \in\Yh_k(u;\theta)\}.
\end{aligned}
\end{equation}
The set $\widehat{\mathcal X}_k(u;\theta,D_v)$ contains
states consistent with the estimated output set and admissible
measurement noise. We consider output and noise sets which can be written as a Cartesian product of sets over sensors, i.e.,  $\Yh_k=\prod\limits_{s \in \mathcal{S}}\Yh_{k,s}$ and $\Vk=\prod\limits_{s \in \mathcal{S}}\mathcal V_{k,s}$, which yields the state set:
\begin{equation}\label{eq:intersection}
\widehat{\mathcal X}_k(u;\theta, D_v)=\bigcap_{s \in \mathcal{S}}
\{x\in\XX:\exists \ v_s \in\mathcal V_{k,s}(D_v),\ C_s x+ v_s\in\Yh_{k,s}(u;\theta)\},
\end{equation}
where the noise vector $v$ in \eqref{eq:preimage} denotes the compact form $v=\operatorname{col}(v_1,\ldots,v_S)\in
\Vk$.  The Cartesian product structure of these sets does not
require statistical independence across sensors. Let $\phi_k:\mathcal U^k\to\mathbb R^{q_k}$
and
$\psi_k:\mathcal U^k\to\mathbb R_{\geq0}^{r_k}$
be feature maps used to parameterize the center and
half-width of the estimated output boxes as follows:
\begin{equation}\label{eq:boxes}
\left\{
\begin{aligned}
c_k(u;\theta)&:=A_k\phi_k(u_{0:k-1})+b_k,  \ \ \ \ \ \ a_k(u;\theta):=a_k^0+\Gamma_k\psi_k(u_{0:k-1}),\\
\Yh_k(u;\theta)&:=\{z\in\R^{n_y}:
|z_i-c_{k,i}(u;\theta)|\le a_{k,i}(u;\theta),
\ i=1,\ldots,n_y\}.
\end{aligned}
\right.
\end{equation}
Here, $c_{k,i}(u;\theta)$ and $a_{k,i}(u;\theta)$
denote the $i$-th components of $c_k(u;\theta)$ and $a_k(u;\theta)$, respectively.
The feature maps $\phi_k$ and $\psi_k$
are assumed to be Borel measurable.
The parameter vector $\theta$ collects
$(A_k,b_k,a^0_{k},\Gamma_k)_{k\in\mathcal T}$,
where
$A_k\in\mathbb R^{n_y\times q_k}$,
$b_k,a^0_{k}\in\mathbb R^{n_y}$,
and
$\Gamma_k\in\mathbb R^{n_y\times r_k}$.
The admissible parameter domain is $\Theta:=
\left\{
\theta:
 a^{0}_k\geq 0,\;
\Gamma_k\geq0,\ 
\forall k\in\mathcal T
\right\}.$
Since $\psi_k\geq0$, the half-widths
remain nonnegative. These linear constraints
define a deterministic convex parameter
domain. For inputs $u(t)\in\R^{n_u}$, a simple choice for the feature maps is $\phi_k(u_{0:k-1})
=
\operatorname{col}(u(0),\ldots,u(k-1))$, and
$\psi_k(u_{0:k-1})
=
\operatorname{col}(|u(0)|,\ldots,|u(k-1)|)$.  For $k\geq 1$, the feature dimensions are
$q_k=r_k=kn_u$.  At $k=0$, we set $q_0=r_0=1$ and
$\phi_0(\emptyset)=\psi_0(\emptyset)=0$,
so that $c_0=b_0$ and $a_0=a^0_{0}$.

We now characterize the challenge  induced
by the sensing map. Specifically, for a state $x\in\mathcal X$, let us define the measurement equivalent  class as
$\mathcal F_C(x):=\{z\in\mathcal X:Cz=Cx\}
=\{x+q:q\in\mathbb R^{n_x},\,Cq=0,\,x+q\in\mathcal X\}$.
Note that states $\tilde{z} \in \mathcal F_C(x)$ produce identical noiseless measurements and
therefore cannot be distinguished using the sensing map alone.
The following result characterizes the
consequences of this challenge for set-based
state estimation.
\begin{proposition}
\label{prop:indistinguishability}
Let $\widehat{\mathcal X}_k(u;\theta,D_v)$ be the
estimated state set. Then, for every
$x\in\widehat{\mathcal X}_k(u;\theta,D_v)$, it holds that $\mathcal F_C(x)
\subseteq
\widehat{\mathcal X}_k(u;\theta,D_v)$. If, in particular, it holds that  $\mathcal X=\mathbb R^{n_x}$, 
$\ker C\neq\{0\}$, and
$\widehat{\mathcal X}_k(u;\theta,D_v)\neq\varnothing$,
then $\operatorname{diam}
\bigl(\widehat{\mathcal X}_k(u;\theta,D_v)\bigr)
=+\infty$. \hfill $\square$
\end{proposition}
\begin{proof}
Let $x\in\widehat{\mathcal X}_k(u;\theta,D_v)$, which implies that there exists
$v\in\mathcal V_k(D_v)$ such that $Cx+v\in\widehat{\mathcal Y}_k(u;\theta)$.
For any $z\in\mathcal F_C(x)$, it holds that
$Cz=Cx$, and thus $Cz+v
\in\widehat{\mathcal Y}_k(u;\theta)$, which in turn implies that
$z\in\widehat{\mathcal X}_k(u;\theta,D_v)$.
If $\mathcal X=\mathbb R^{n_x}$
and $\ker C\neq\{0\}$, let us choose $q\in\ker C\setminus\{0\}$. Then for every $t\in\mathbb R$, it holds that  $C(x+tq)=Cx$,
so that $x+tq\in\widehat{\mathcal X}_k(u;\theta,D_v)$, i.e., $\widehat{\mathcal X}_k(u;\theta,D_v)$ is unbounded.
\end{proof}

Due to Proposition \ref{prop:indistinguishability}, we impose the following assumption: 
\begin{assumption} \label{assum:rank}
The rank of $C$ is equal to the dimension of the state vector, i.e., $\rank(C)=n_x$. \hfill $\square$
\end{assumption}

Assumption~\ref{assum:rank}
ensures that distinct states produce distinct
noiseless outputs. Combined with
bounded output and measurement noise sets it guarantees finite state set diameter bounds without knowing a model of the state dynamics in \eqref{eq:state_dynamics}. Let $C^\dagger=(C^\top C)^{-1}C^\top$, such that $C^\dagger C=I$ holds. For a nonempty bounded $\Vk$, let us define $\omega_i:=\|C^\dagger e_i\|_2$ and $b_{v,k}:=\sup\limits_{v_1,v_2\in\Vk}\|C^\dagger(v_1-v_2)\|_2$.
\begin{proposition}\label{prop:surrogate}
For any input $u$ and
$\theta\in\Theta$ such that
$\widehat{\mathcal X}_k(u;\theta,D_v)$
is nonempty, it holds that
$\operatorname{diam}\bigl(
\widehat{\mathcal X}_k(u;\theta,D_v)
\bigr)
\leq
2\sum\limits_{i=1}^{n_y}
\omega_i a_{k,i}(u;\theta)
+b_{v,k}(D_v).$ \hfill $\square$
\end{proposition}
\begin{proof}
Let $x_1,x_2\in
\widehat{\mathcal X}_k(u;\theta,D_v)$.
By construction, there exist
$z_1,z_2\in\widehat{\mathcal Y}_k(u;\theta)$
and $v_1,v_2\in\mathcal V_k(D_v)$ such that
$Cx_\ell=z_\ell-v_\ell$ for $\ell=1,2$.
Since $C^\dagger C=I$,
\[
\begin{aligned}
\|x_1-x_2\|_2 \leq
\sum_{i=1}^{n_y}
\|C^\dagger e_i\|_2
|z_{1,i}-z_{2,i}| +
\|C^\dagger(v_1-v_2)\|_2 \leq
2\sum_{i=1}^{n_y}
\omega_i a_{k,i}(u;\theta)
+b_{v,k}(D_v).
\end{aligned}
\]
The last inequality follows from
$|z_{1,i}-z_{2,i}|\leq
2a_{k,i}(u;\theta)$.
Taking the supremum over all pairs
$x_1,x_2$ concludes the proof.
\end{proof}
Proposition~\ref{prop:surrogate}
establishes an upper bound on the diameter
of the state set as a function of
the  half-widths of the estimated output set and the
geometry of the sensor matrix $C$. The coefficients
$\omega_i=\|C^\dagger e_i\|_2$
quantify the contribution of uncertainty
in each output coordinate to the state set diameter. We have now established the core ingredients for
quantifying the estimated state sets under unknown
dynamics, which we will use to provide probabilistic
guarantees that these sets contain the true state
of the system in \eqref{eq:system}.
\section{Set-based State Estimation
via Scenario Optimization with Relaxation}
\label{sec:certificate} Proposition~\ref{prop:surrogate} establishes
an upper bound on the diameter of the
estimated state sets in terms of the size of the output set and the sensor
matrix $C$. 
To obtain a computationally tractable
learning formulation, we introduce
the diameter surrogate $d_k(u;\theta)
:=
2\sum_{i=1}^{n_y}
\omega_i a_{k,i}(u;\theta)$. 
Since $b_{v,k}(D_v)$ is independent
of the learning parameters $\theta$,
we minimize $d_k(u;\theta)$ to reduce
the parameter-dependent upper bound
on the estimated state set diameter.
Let $y_i(k)\in\mathbb R$ denote the $i$-th component
of the concatenated output $y(k)\in\mathbb R^{n_y}$
for $i=1,\ldots,n_y$.
To quantify output set violations
at the trajectory level, we define
the uncertain function $g(\zeta,\theta):=
 \max\limits_{\substack{k\in\mathcal T, 1\le i\le n_y}}
 \bigl\{|y_i(k)-c_{k,i}(u;\theta)|-a_{k,i}(u;\theta)\bigr\}$.
Note that the constraint $g(\zeta,\theta)\leq0$
is equivalent to output inclusion at
every time step. To accommodate trajectories that may violate
the estimated output sets, we introduce
a nonnegative relaxation variable $\xi_j$ for each
training trajectory.
Denoting the sample index set $\mathcal N:=\{1,\ldots,N\}$, for a fixed penalty
$\rho>0$, we solve the scenario program
\begin{equation}
\left\{
\begin{aligned}
\min_{\theta\in\Theta,\,\xi\ge0}\quad &
\sum_{k\in\mathcal T}
\max_{j\in\mathcal N}d_k(u^{(j)};\theta)
+\rho\sum_{j \in \mathcal{N}}\xi_j,\\
\text{subject to}\quad &
g(\zeta^{(j)},\theta)\le\xi_j,
\ j\in\mathcal N. \nonumber 
\end{aligned}
\right.
\end{equation}

The cost balances the worst-case diameter surrogate over
the samples $\zeta^{(j)}$, with $j\in\mathcal N$, at each time $k \in \mathcal{T}$
against output relaxation for each trajectory. Introducing the epigraphic variables  $\tau=(\tau_k)_{k\in\mathcal T} \in \mathbb R_{\geq0}^{T+1}$ we obtain the following epigraphic reformulation which is a linear program (LP) due to the geometry of the output set.
\begin{equation}
\left\{
\begin{aligned}
\min_{\theta\in\Theta,\,\tau\ge0,\,\xi\ge0}\quad &
\sum_{k\in\mathcal T}\tau_k
+\rho\sum_{j \in \mathcal{N}}\xi_j\\
\text{subject to}\quad &
|y_i^{(j)}(k)-c_{k,i}(u^{(j)};\theta)|
\leq a_{k,i}(u^{(j)};\theta)+\xi_j,\\
&d_k(u^{(j)};\theta)\le\tau_k, \  j\in\mathcal N,\ k\in\mathcal T, \  i=1,\ldots,n_y.
\end{aligned}
\right.
\label{eq:relaxed_scenario_reformulation}
\end{equation}
Note that the epigraph constraints on $d_k$
are hard constraints, while constraints corresponding to the inclusion of output trajectories are
relaxed through the variables $\xi_j$.
For fixed training inputs, \eqref{eq:relaxed_scenario_reformulation} is linear in the
decision variables and feasible for every finite
multi-sample $D_N$.
To obtain a unique optimizer, we select
the minimum-norm solution among all
optimizers of
\eqref{eq:relaxed_scenario_reformulation},
by minimizing $\|(\theta,\tau,\xi)\|_2^2$
over the optimal set of
\eqref{eq:relaxed_scenario_reformulation}. The corresponding optimal epigraph cost and relaxation variables then satisfy $\tau_{N,k}^\star= \max_{j\in\mathcal N}d_k(u^{(j)};\theta_N^\star)$, $\xi_{N,j}^\star=\max{(g(\zeta^{(j)},\theta_N^\star), 0)}$, respectively. For every $k\in\mathcal T$, the estimated state set
is obtained by evaluating (2) for
$\theta=\theta_N^\star$, namely,
\begin{equation}
\begin{aligned}
\widehat{\mathcal X}_k(u;\theta_N^\star,D_v)
&:=
\left\{
x\in\mathcal X:
\exists v\in\mathcal V_k(D_v),\;
Cx+v\in
\widehat{\mathcal Y}_k(u;\theta_N^\star)
\right\}.
\end{aligned}
\end{equation}
In particular, for box-shaped noise sets
$\mathcal V_k(D_v)=
\{v:|v|\leq\nu_k(D_v)\}$,
this reduces to
$\widehat{\mathcal X}_k(u;\theta_N^\star,D_v)
=
\left\{
x\in\mathcal X:
|Cx-c_k(u;\theta_N^\star)|
\leq a_k(u;\theta_N^\star)+\nu_k(D_v)
\right\}. $
Here, absolute values and inequalities
are understood componentwise.

\section{Probabilistic Certificates for Set-based State Estimation}
\label{pac:sec:certificates}

We first seek distribution-free probabilistic guarantees
that a new trajectory satisfies both output inclusion
and the learned diameter thresholds.
We achieve this using the scenario optimization
framework of \cite{GarattiCampi2025}.
First, we show that their Theorem~4 applies to
problem \eqref{eq:relaxed_scenario_reformulation}
by verifying the consistency conditions in
their Property~1. Let $\Xi:=\UU^T\times\R^{(T+1)n_y}$
be the space of input-output trajectories and define
the decision domain
$\mathcal Z:=\Theta\times\R_{\geq0}^{T+1}
\times\mathbb Z_{\geq0}$.

\begin{definition} (Adapted from \cite{GarattiCampi2025})
\label{def:decision-appropriateness}
For a fixed penalty $\rho>0$, consider for each $m\geq1$
an ordered multi-sample
$\zeta_m=(\zeta^{(1)},\ldots,\zeta^{(m)})
\in\Xi^m$. Let
$(\theta_m^\star,\tau_m^\star,\xi_m^\star)$
be the minimum-norm optimizer of
\eqref{eq:relaxed_scenario_reformulation}
when only these $m$ samples are used.
Let
$r_m:=\sum_{j=1}^{m}
\mathbf 1_{\{g(\zeta^{(j)},\theta_m^\star)>0\}}$
denote the number of sampled trajectories that violate
output inclusion.
The decision map $M_m:\Xi^m\to\mathcal Z$ is defined by
\begin{equation}\label{pac:eq:decision}
M_m(\zeta_m):=
(\theta_m^\star,\tau_m^\star,r_m).
\end{equation}
For the empty multi-sample, define
$M_0(\emptyset):=(0,0,0)$,
corresponding to the minimum-norm solution of
\eqref{eq:relaxed_scenario_reformulation}
with no sample constraints and an empty slack vector.
 For each trajectory $\zeta=(u,y)\in\Xi$, the appropriateness set is given by
\begin{equation}\label{pac:eq:appropriateness}
\mathcal A_\zeta
:=
\left\{
(\theta,\tau,r)\in\mathcal Z:
g(\zeta,\theta)\leq0, \
d_k(u;\theta)\leq\tau_k
 \text{ for every }k\in\T
\right\}.
\end{equation}
The coordinate $r$ represents the number of training trajectories that violate output inclusion. A decision $z\in\mathcal Z$ is appropriate for $\zeta$
if $z\in\mathcal A_\zeta$. \hfill $\square$
\end{definition}
Evaluating $M_m(\zeta_m)$ means solving (5)
on these $m$ trajectories and returning the learned
parameters, diameter surrogate thresholds, and number
of training trajectories which violate the output set.
A decision belongs to $\mathcal A_\zeta$ if output
inclusion holds without relaxation and all diameter
surrogate thresholds are satisfied.  For a subset $I\subseteq\mathcal N$, let $D_I$
be the sublist of $D_N$ indexed by $I$,
in increasing order, and write
$z_I:=(\theta_I^\star,\tau_I^\star,r_I)
:=M_{|I|}(D_I)$. Thus, $r_I$ counts the trajectories in $D_I$ that violate
output inclusion, and
$z_N:=(\theta_N^\star,\tau_N^\star,r_N)=M_N(D_N)$.

\begin{lemma}\label{lem:consistency}
For every $m\geq0$, $n\geq1$, and
$\zeta_{m+n}\in\Xi^{m+n}$, write
$z_m:=M_m(\zeta_m)$ and
$z_{m+n}:=M_{m+n}(\zeta_{m+n})$,
where $\zeta_m$ contains the first $m$ trajectories.
The maps satisfy the three conditions of
\cite[Property~1]{GarattiCampi2025}:  (i)
$M_m(\zeta^{(\pi(1))},\ldots,\zeta^{(\pi(m))})=z_m$
for every permutation $\pi$ of $\{1,\ldots,m\}$,
(ii)
$z_m\in\bigcap_{j=m+1}^{m+n}\mathcal A_{\zeta^{(j)}}$ implies that $z_{m+n}=z_m$ and
(iii)
$z_m\notin\bigcap_{j=m+1}^{m+n}\mathcal A_{\zeta^{(j)}}$
 implies that $z_{m+n}\neq z_m$.
\end{lemma}
\begin{proof}
 (i) Permuting samples
does not change the value of the map leaving the  decision unchanged.
(ii) Let $J_\ell^\star$ denote the optimal value of \eqref{eq:relaxed_scenario_reformulation}
with $\ell$ trajectories. Since every added trajectory
is appropriate,
$(\theta_m^\star,\tau_m^\star,\xi_m^\star,0_n)$
is feasible with $m+n$ trajectories, i.e., 
$J_{m+n}^\star\leq J_m^\star$.
Conversely, any feasible $(\theta,\tau,\xi)$ for these
$m+n$ trajectories yields a feasible
$(\theta,\tau,\xi_{1:m})$ for the original problem.
Its objective is therefore at least
$J_m^\star+\rho\sum_{j=m+1}^{m+n}\xi_j\geq J_m^\star$, which implies that $J_{m+n}^\star=J_m^\star$.
Since $\rho>0$, every optimizer has zero added slacks
and restricts to an original optimizer. Appending zeros preserves the norm, so the minimum-norm
rule selects
$(\theta_m^\star,\tau_m^\star,\xi_m^\star,0_n)$.
The added trajectories contribute no violations, i.e.,  $z_{m+n}=z_m$.
(iii) If an added diameter constraint fails at
$(\theta_m^\star,\tau_m^\star)$, this pair becomes
infeasible and must change,
otherwise, some added trajectory has
$g(\zeta^{(j)},\theta_m^\star)>0$.
If the selected pair stays unchanged, then
$r_{m+n}=r_m+\sum_{j=m+1}^{m+n}
\mathbf 1_{\{g(\zeta^{(j)},\theta_m^\star)>0\}}>r_m$.
Thus it holds that $z_{m+n}\neq z_m$.
\end{proof}
\subsection{Output coverage and size certification}
\label{pac:subsec:output-certificate}

Let us define the complexity  $s_N:=\min\left\{|I|:
 \begin{array}{l}
 I_+\subseteq I\subseteq\mathcal N,
 \theta_I^\star=\theta_N^\star, \
 \tau_I^\star=\tau_N^\star
 \end{array}\right\}$,
 where $I_+:=\{j\in\mathcal N:g(\zeta^{(j)},\theta_N^\star)>0\}$. 
The requirement $I^+\subseteq I$ ensures that
whenever the selected parameters and diameter
thresholds are reconstructed, the count of violations from the original scenario program with all trajectories in place is also preserved. Thus, the subset $I$ reconstructs the 
decision $z_N=M_N(D_N)$, obtained when all samples are considered. For a new input/output trajectory $\zeta=(u,y)$ from $\mathbb P_0$,
we define the following probabilities:
\begin{align}
&\tilde{\mathbb V}(D_N)
:=\mathbb P_0\left\{
\zeta:
g(\zeta,\theta_N^\star)>0
\ \text{or}\
\exists k\in\mathcal T:
d_k(u;\theta_N^\star)>\tau_{N,k}^\star
\right\}, \nonumber \\
&\mathbb V_y(D_N)
:=\mathbb P_0\left\{
\zeta:g(\zeta,\theta_N^\star)>0
\right\}. \nonumber 
\end{align}
The function $\mathbb V_y(D_N)$ is the probability
that for some time $k \in \mathcal{T}$, a new yet unseen output measurement does not belong in the estimated output set $\widehat{\mathcal Y}_k(u;\theta_N^\star)$. The function
$\tilde{\mathbb V}(D_N)$ is the probability
that either output inclusion fails or the
diameter surrogate exceeds the 
threshold obtained by the scenario program in \eqref{eq:relaxed_scenario_reformulation}, at some time step for a new yet unseen realization $\zeta=(u,y) \sim \mathbb{P}_0$.
Consequently, it holds that
$\mathbb V_y(D_N)\leq
\tilde{\mathbb V}(D_N)$. For $\beta\in(0,1)$ and $q\in\{0,\ldots,N-1\}$,
let $t_q$ be the unique root in $(0,1)$ of the
following polynomial in $t$:
\begin{equation}
 \binom{N}{q}t^{N-q}
 -\frac{\beta}{N}\sum_{m=q}^{N-1}\binom{m}{q}t^{m-q}=0.
 \label{pac:eq:epsilon-polynomial}
\end{equation}
Let $\mathbb P_0^N:=\prod_{i \in \mathcal{N}}\mathbb P_0$
denote the distribution of $D_N$.
The following result then holds:

\begin{proposition}\label{pac:thm:output}
Suppose Assumption~\ref{ass:law} holds and the maps
$M_m:\Xi^m\to\mathcal Z$ are Borel measurable
for every $m\geq0$. Let $s_N$ be the complexity
defined above. For $\beta\in(0,1)$, define
$\overline\epsilon_N(q,\beta):=1-t_q$ for
$q\in\{0,\ldots,N-1\}$, where $t_q\in(0,1)$
is the unique root of (9), and set
$\overline\epsilon_N(N,\beta):=1$.
Then, it holds that
\begin{equation}\label{pac:eq:augmented-certificate}
\mathbb{P}_0^N\left\{
D_N:\tilde{\mathbb{V}}(D_N)
\leq\overline\epsilon_N(s_N,\beta)
\right\}
\geq 1-\beta.
\end{equation}
\end{proposition}

\begin{proof}
By Lemma~\ref{lem:consistency}, the maps $M_m$
satisfy \cite[Property~1]{GarattiCampi2025}.
The result follows from
\cite[Theorem~4]{GarattiCampi2025}.
\end{proof}
Since $\epsilon_N(q,\beta)$ is nondecreasing in $q$
\cite[Section~2.1]{GarattiCampi2025}, any measurable,
integer-valued upper bound $\widehat s_N$ satisfying
$s_N\leq\widehat s_N\leq N$ yields a valid certificate. Such a bound can be obtained by
sequentially removing trajectories outside
$I_+$, retaining a removal only if resolving
\eqref{eq:relaxed_scenario_reformulation}
with the same penalty and tie-break rule
preserves $(\theta_N^\star,\tau_N^\star)$.
The cardinality of the remaining subset
then provides $\widehat{s}_N$.
\subsection{Probabilistic Guarantees for State Inclusion and Set Tightness}
\label{pac:subsec:state-certificate}

We remind the reader that the variable $E=(x_0,w_{0:T-1},v_{0:T})
\sim\mathbb P_E$ collects the uncertainties
affecting the system in \eqref{eq:state_dynamics}-\eqref{eq:output_dynamics}, whose corresponding
input-output probability distribution is $\mathbb{P}_0$. Each realization of $E$ induces a state, process and measurement noise trajectory under the considered input.  For notational simplicity, with a slight abuse of notation, we write
$x(k)=x(k;E)$, $y(k)=y(k;E)$, and
$u=u(E)$, with the dependence on the
uncertainty realization $E$ to be implied.
Let $\mathbb P_D$ denote the joint distribution of the multi-sample
$D=(D_N,D_v)$. For any such multi-sample, we define
the state trajectory exclusion probability
\begin{equation}
\mathbb{V}_x(D):=\mathbb P_E\left\{
E:\exists k\in\mathcal T,\;
x(k)\notin
\Xh_k(u;\theta_N^\star,D_v)
\right\}. \nonumber 
\end{equation}
The function $\mathbb V_x(D)$ is the probability that a new
state trajectory violates the estimated sets at some time
$k\in\mathcal T$. The following result provides a bound on this probability.

\begin{proposition}\label{pac:prop:transfer}
For any $D$, it holds that
$\mathbb{V}_x(D)
\leq \mathbb{V}_y(D_N)+\mathbb{V}_v(D_v)
\leq \widetilde{\mathbb{V}}(D_N)+\mathbb{V}_v(D_v)$.
\hfill $\square$
\end{proposition}
\begin{proof}
Note that if
$y(k)\in\widehat{\mathcal Y}_k(u;\theta_N^\star)$
and $v(k)\in\mathcal V_k(D_v)$,
then $x(k)\in\widehat{\mathcal X}_k(u;\theta_N^\star,D_v)$.
Hence, by the subadditivity property of the union of two events,
\begin{align*}
\mathbb{V}_x(D)
&\leq
\mathbb P_E\left\{
E:\exists k\in\mathcal T,\;
y(k)\notin\widehat{\mathcal Y}_k(u;\theta_N^\star)
\right\} +
\mathbb P_E\left\{
E:\exists k\in\mathcal T,\;
v(k)\notin\mathcal V_k(D_v)
\right\}\\
&=\mathbb{V}_y(D_N)+\mathbb{V}_v(D_v) \leq \widetilde{\mathbb{V}}(D_N)+\mathbb{V}_v(D_v).
\end{align*}
The last inequality follows from
$\mathbb{V}_y(D_N)\leq\widetilde{\mathbb{V}}(D_N)$,
which concludes the proof.
\end{proof}

For each time step $k\in\mathcal T$, we define the threshold
$B_k(D):=\tau_{N,k}^\star+b_{v,k}(D_v)$.
Let us now define the probability of violation for joint state coverage and tightness of the state set as follows:
\begin{equation}
\begin{aligned}
\mathbb V_{\mathrm{cs}}(D)
:=\mathbb P_E\Bigl\{
E:\exists k\in\mathcal T:\;
&x(k)\notin\widehat{\mathcal X}_k(u;\theta_N^\star,D_v) \text{ or }
\operatorname{diam}
\bigl(\widehat{\mathcal X}_k(u;\theta_N^\star,D_v)\bigr)>B_k(D)
\Bigr\}. \nonumber 
\end{aligned}
\end{equation}
Unlike $\mathbb V_x(D)$, which measures only state trajectory exclusion, $\mathbb V_{\mathrm{cs}}(D)$ also accounts for violations of the obtained
diameter threshold. For example, a set may
contain the true state while exceeding the
allowed diameter, i.e,  state inclusion guarantees alone do not ensure satisfaction of diameter bounds.

\begin{theorem}\label{pac:thm:state}
Under Assumptions~1-5, let  the maps
$M_m:\Xi^m\to\mathcal Z$ be Borel measurable
for every $m\geq0$. Let
$\overline\epsilon_N$ be defined in
Proposition~\ref{pac:thm:output}, and
$\epsilon_v(D_v,\beta_v)$ as in Assumption~4.
Then, for any $\beta_y\in(0,1-\beta_v)$ and complexity bound
$\widehat s_N$ satisfying
$s_N\leq\widehat s_N\leq N$, it holds that
\begin{equation}\label{pac:eq:state-certificate}
\mathbb P_D\left\{
D:\mathbb V_{\mathrm{cs}}(D)
\leq
\overline{\epsilon}_N(\widehat{s}_N,\beta_y)
+\epsilon_v(D_v, \beta_v)
\right\}
\geq 1-\beta_y-\beta_v.
\end{equation}
\end{theorem}

\begin{proof}
Fix $D$ and consider a realization $E$ for which
output and noise inclusion and the
diameter surrogate threshold holds for all
$k\in\mathcal T$.
By construction, we then have that 
$x(k)\in\widehat{\mathcal X}_k(u;\theta_N^\star,D_v)$.
Moreover, by Proposition~\ref{prop:surrogate},  it holds that $\operatorname{diam}\bigl(
\widehat{\mathcal X}_k(u;\theta_N^\star,D_v)
\bigr) \leq d_k(u;\theta_N^\star)+b_{v,k}(D_v) \leq \tau_{N,k}^\star+b_{v,k}(D_v)
=B_k(D).$
Consequently, we have that $\mathbb V_{\mathrm{cs}}(D)
\leq \tilde{\mathbb V}(D_N)
+\mathbb V_v(D_v).$
The result follows from
\eqref{pac:eq:augmented-certificate},
Assumption~\ref{assum:noise_set},
and the subadditivity property.
\end{proof}
Theorem \ref{pac:thm:state} thus provides a finite-sample probabilistic certificate
for simultaneous state inclusion and satisfaction of the learned
diameter bounds throughout the horizon.

\section{Numerical Simulations}\label{sec:numerical}

\begin{figure}[t]
\centering
\includegraphics[width=\linewidth]{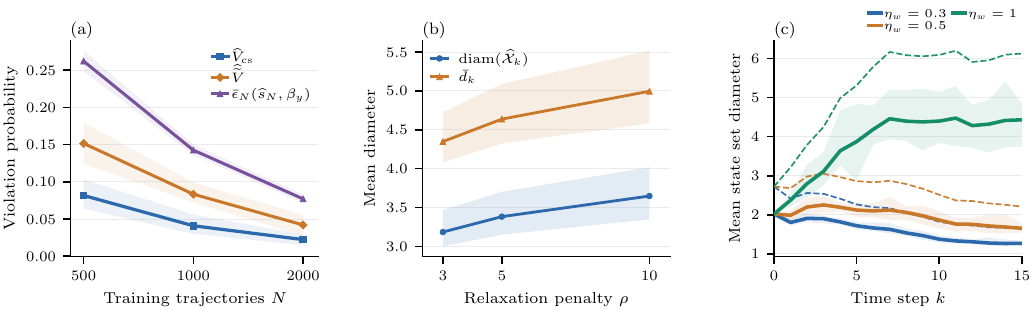}
\caption{\footnotesize (a) Empirical violations $\widehat{ V}_{cs}$,
$\widehat{\widetilde{ V}}$ and bound
$\overline\epsilon_N(\widehat s_N,\beta_y)$
averaged over 15 repetitions for $N \in \{500, 1000, 2000\}$, $\rho=5$ and  $\eta_w=1$. Shaded regions indicate the min/max range across the
15 independent repetitions. (b) State set diameter and bound $\bar d_k:=d_k+b_v$,
both averaged over time steps, test trajectories,
and 15 independent repetitions, for
$\rho\in\{3,5,10\}$, $N=1000$, and $\eta_w=1$. (c)  Evolution of the mean actual and surrogate state set
diameters for $\eta_w\in\{0.3,0.5,1\}$, with $N=2000$
and $\rho=5$. Solid and dashed lines represent the
actual and surrogate diameters, respectively, while
shaded regions indicate the min--max range across
$15$ independent repetitions.
}
\label{fig:numerical-three-panel}
\end{figure} 

We consider a two-dimensional unknown nonlinear system with state dynamics
$x_{k+1}=Ax_k-0.9B_\phi\tanh(x_{1,k})+B_u u_k+w_k$, and output dynamics 
$y_k=Cx_k+v_k$, with
$A=\bigl[\begin{smallmatrix}0.95&0.10\\-0.20&0.85\end{smallmatrix}\bigr]$,
$B_\phi=(0.18,0.06)^\top$, $B_u=(0.15,0.24)^\top$, and
$C=\operatorname{diag}(1,0.14)$.
Initial states are randomly chosen from a  uniform distribution on $[-0.6,0.6]\times[-0.45,0.45]$. The measurement noise is independent and uniform on the support set
$\mathcal V=[-0.10,0.10]\times[-0.014,0.014]$.
The input is given by the noisy output feedback law 
$u_k=\operatorname{sat}_{[-0.65,0.65]}\!
\bigl(0.58u_{k-1}-(0.42,16/7)y_k\bigr)$ and
$u_{-1}=0$. Conditional on the input $u_k$, the
Gaussian process noise components have standard deviations
$\sigma_{w,1}(u_k)=0.037+\eta_w\max\{u_k,0\}$ and
$\sigma_{w,2}(u_k)=0.037+\eta_w\max\{-u_k,0\}$,
where $\eta_w\geq0$ controls their input dependence.  We then consider $\chi_0=0$ and
$\chi_{k+1}=0.8\chi_k+u_k$. The centers of the output
box are affine in $\chi_k$, while their nonnegative
half-widths are affine in
$[\chi_k]_+$ and $[-\chi_k]_+$,
where $[a]_+:=\max\{a,0\}$.
We then solve the scenario program in \eqref{eq:relaxed_scenario_reformulation} with  $d_k(u;\theta)=2a_{k,1}(u;\theta)+\dfrac{2a_{k,2}(u;\theta)}{0.14}$ and $b_v=2\sqrt{0.02}$, which yields $\bar{d}_k=d_k+b_v$ (see Proposition \ref{prop:surrogate}).
For $\mathcal X=\mathbb R^2$, the exact state set diameter, used only for evaluation, is $\operatorname{diam}(\widehat{\mathcal X}_k)
=2\sqrt{(a_{k,1}+0.10)^2+\bigl((a_{k,2}+0.014)/0.14\bigr)^2}$. We use $T=15$, $\beta_y=10^{-4}$,
and the known noise set $\mathcal V$, i.e., $\epsilon_v=\beta_v=0$.
Empirical violations are evaluated on
$M=2500$ independent test trajectories per repetition, over $15$ independent
repetitions. Specifically, the quantities $\widehat{V}_{cs}$ and
$\widehat{\widetilde{V}}$ denote the empirical counterparts of the probabilities $\mathbb V_{cs}$
and $\widetilde{\mathbb V}$, respectively. State observations are used
only for evaluation, not training.   The results in Figure~\ref{fig:numerical-three-panel} indicate
improved generalization to new unseen states with increasing
number of training trajectories, as empirical violations
decrease and probabilistic bounds become less conservative.
For fixed $\rho$, increasing $N$ also increases the weight
on average relaxation relative to the diameter objective. Moreover, the relaxation penalty tunes
the trade-off between state coverage and state set tightness. 
Finally, the decreasing mean state set diameter
under weaker disturbances suggests that the
attenuation of initial uncertainty dominates
the accumulation of process noise, whereas
stronger disturbances prevent such narrowing.

\section{Concluding Remarks}

We developed a data-driven framework for
set-based state estimation from noisy input and
output trajectories without identifying the
unknown dynamics or requiring state measurements.
We derived diameter bounds based on sensing
geometry and formulated a tractable relaxed
linear program. Leveraging the results in scenario optimization,
we established a posteriori probabilistic
guarantees for state inclusion and satisfaction of diameter bounds. Future work will investigate nonlinear
sensing maps 
and extensions to secure set-based state estimation.

\section*{Acknowledgments}
\emph{\small During the preparation of this work, the authors used
ChatGPT to improve readability, assist with grammar
and phrasing, and support code development.
The authors reviewed and edited the text and code
as needed and take full responsibility for the
content of the publication.}

\bibliography{references}
\end{document}